# Feature Selection Approaches for Newborn Birthweight Prediction in Multiple Linear Regression Models

Esther Liu
Pei Xi Lin
Qianqi Wang
Karina Chen Feng

# Abstract

This project is based on the dataset "exposome_NA.RData", which contains a subcohort of 1301 mother-child pairs who were enrolled into the HELIX study during pregnancy. Several health outcomes were measured on the child at birth or at age 6-11 years, taking environmental exposures of interest and other covariates into account.

This report outlines the process of obtaining the best MLR model with optimal predictive power. We first obtain three candidate models we obtained from the forward select, backward elimination and stepwise selection, and select the optimal model using various comparison schemes including AIC, Adjusted R^2 and cross validation for 8000 repetitions.

The report ended with some additional findings revealed by the selected model, along with restrictions of the method we use in the model selection process.

# Objectives

The objective of the analysis is to provide the best predictive model for birthweight using multiple linear regression models. This report outlines the process of obtaining the three candidate models obtaining the best MLR model with optimal predictive power by considering the model selection process using cross validation.

# Exploratory Data Analysis

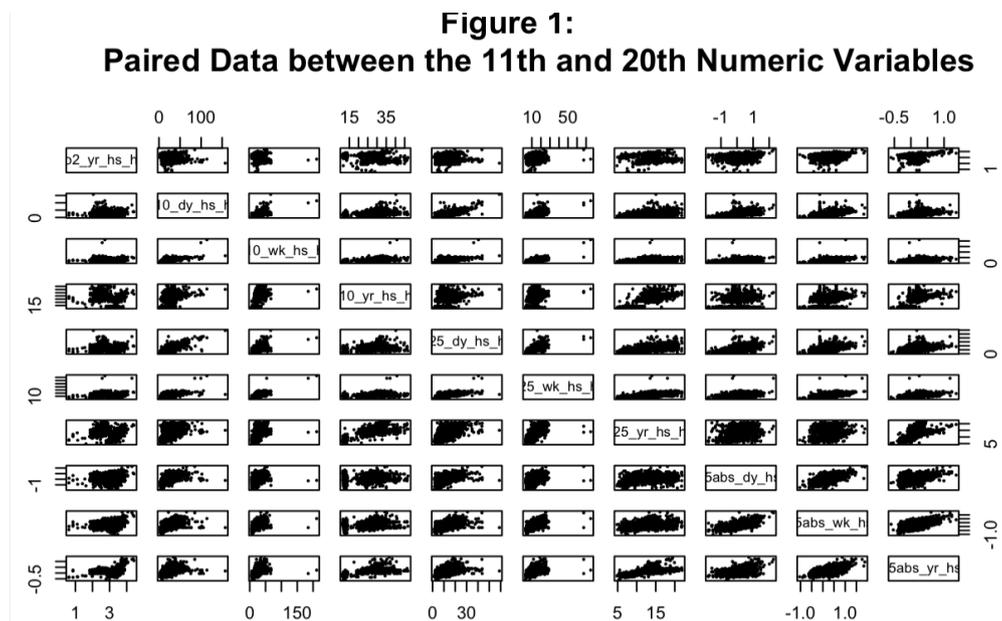

Figure 1:
Paired Data between the 11th and 20th Numeric Variables

Multicollinearity describes the scenario where multiple parameters are highly correlated to one another. It is difficult to determine the actual effect of each parameter on the response variable with the existence of multicollinearity. As a result, the inclusion and the estimation of variables towards the model can be inaccurate. As we see above, a scatterplot matrix can point out the pairs of variables that are highly correlated, i.e. the 16th and 17th numeric variables. One method we considered to detect the existence of multicollinearity in our model is by the variance inflation factor, or VIF. It is defined as a measurement over the inflation, caused by multicollinearity, of the standard error of the estimate of the coefficient. In most scenarios, there will be multicollinearity to some extent according to VIF. By experience, we say that a VIF that is greater than 10 signifies a problematic multicollinearity. Our solution is thus to remove the redundant variable with high VIF from the model, until the value of VIF is controlled below 10 as desired.

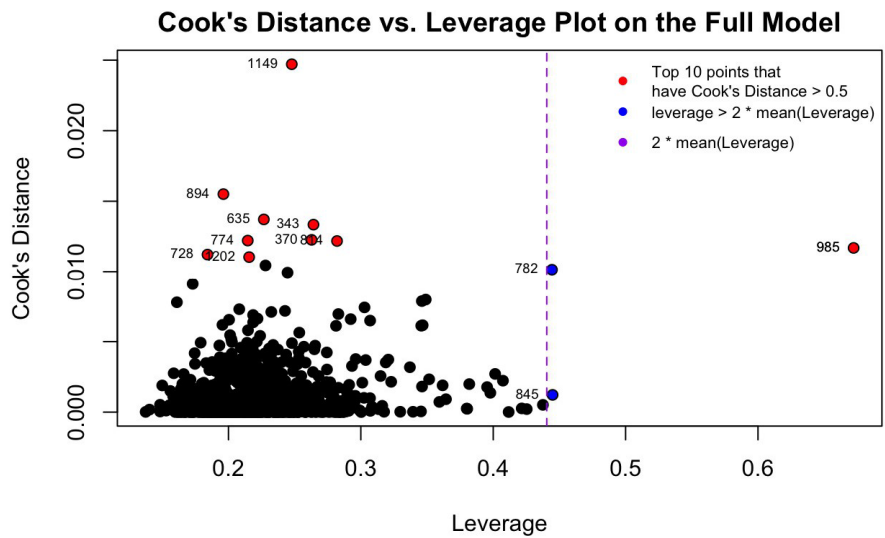

**Figure 2. Cook's Distance vs Leverage Scatterplot for the Full Model lm(e3_bw~ . , data=predict_data)**

Also, by looking at the influence and leverage for each observation in the above scatterplot (where 10 observations with the highest Cook's distance are labelled red), we noticed that there are some points with a leverage greater than or equal to 2 times the average leverage. Specifically, the leverage for observation 985 is significantly higher than the rest of the observations, and it is also in the top 10 points with highest Cook's distance. This indicates that the 985th observation is an outlier. Instead of excluding it, we will examine how the results will look like with and without the 985th observation, and see if our conclusions are substantively different for good practice.

# Methods

This investigation considered all variables for model selection with exceptions for the phenotypes from which only birthweight was considered. This is because phenotypeNA describes several health outcomes measured and only birthweight, as a result variable, is concerned. Therefore, the decision is to not include other result variables into the research, and only use variates from the covariates and exposome tables.

Another notable decision made was the transformation of some variables from numeric types to factor types. During the data preview and cleansing process, some flag-like variables were noticed to have only binary values, i.e. ones and zeros. These variables are in fact discrete numbers in the data, however, it is possible that R mistakenly treats such variables as continuous covariates to maintain linearity. Therefore, converting such variables to factor type allows R to assign appropriate values and recognize these categorical variables. Also, necessary numbers of indicator variables are assigned after this conversion. For instance, for a categorical variable with n levels, R will include n-1 indicator variables in the model.

We then perform iterative elimination of numerical variates with high variance influence factor (VIF) to reduce potential multicollinearity. We set the $v_{*} = 10$ as the threshold for a large VIF by convention (Vittinghoff et al, 2011). This eliminates 7 variables and reduces the number of variables from 229 to 222 (number of columns being reduced from 230 to 223, since the dataset includes the ID column).

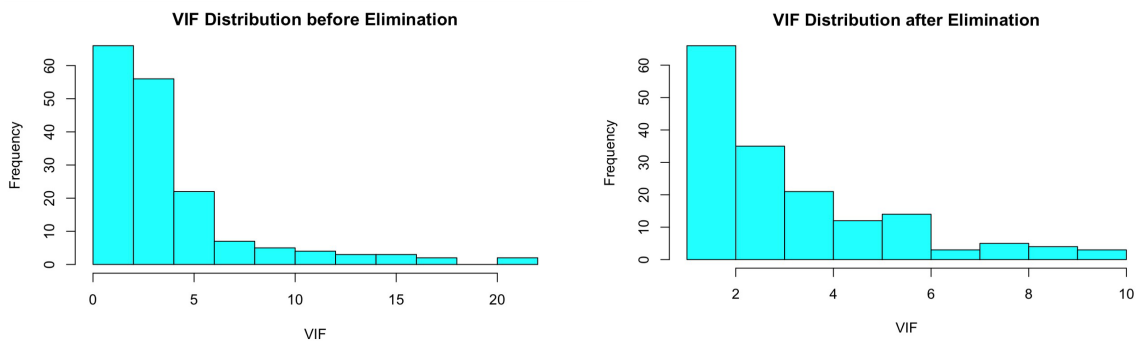

**Figure 3.**
**Distributions of VIFs for Numeric Observations Before and After Elimination**

The histograms show the distribution of the VIFs of numeric data before and after the elimination process. Since the number of variables is still large after this elimination, we will then consider three different model selection processes: Forward selection, backward elimination and stepwise selection. We use AIC as the criteria for all three model selection methods, by setting the degree of freedom used for the penalty k to 2 in the "step" function.

In **forward selection**, the model starts with an empty model and adds in variables individually. The one variable that gives the best improvement to the current model is added in each step. To begin with, the model is formed with only one intercept. Various variables are tested to verify their relevance to the model, where the very best variable is determined according to one fixed criteria. It is such a variable that is added to the model each time. As the model continues to improve along with such testing processes, variables are added in the model one by one. At the point when it is impossible for the model to better improve, the model is considered optimised and the process stops.

In **backward selection**, the model starts with a model including all variables and removes variables one by one. The one variable that is the least significant in terms of a certain criteria is removed each time. To begin with, the regression is formed onto all variables. They are tested to determine their significance to the model, where the least useful one is what to be removed. As the model continues to improve closer to the desired state, nonsignificant variables are removed from the model one by one. At the point when the pre-specified criteria is satisfied, the removing process may stop and the model is considered to be at an ideal state.

In **stepwise selection**, the model starts with only an intercept, excluding all covariates at first. Using the same idea as forward selection, we add one variate to the model. After adding a covariate to the model, we check our new model to see if any covariates should be removed using the same steps in backward selection. We repeat the steps of adding new covariates, and checking whether any covariates need to be removed until all the covariates have been covered. Stepwise selection can be viewed as performing forward and backwards at the same time on all the covariates.

## Model Diagnostics

We have selected three candidate models, and we will compare some fundamentals of them. We will begin with the following metrics:

```
                        Forward         Backward        Stepwise
sum-of-squared PRESS  263388165.3889  263544033.0570  263586592.71
AIC                      19234.3256      19234.2113     19234.19
Adjusted R^2                 0.2435          0.2579         0.26
sum-of-squared DFFITS       71.6536        100.2147       104.65
Number of Predictors        67.0000         93.0000        97.00
```

We noticed that the Stepwise selection model is preferred with respect to the adjusted Coefficient of determination (adjusted $R^2$) and Akaike Information Criterion (AIC), while the Forward selection model with respect to the sum-of-squared PRESS residuals and the sum-of-squared DFFITS. We also noticed that there is a large difference in sum of squared DFFITS of Forward selection model and the sum of squared DFFITS of the Backward and Stepwise model. This suggests that the Backward and Stepwise model is being affected by high leverage observations more severely than the Forward selection model.
Hence, we investigate the high-leverage points and their influence towards all three candidate models. The below figure shows the Cook's Distance vs Leverage for all three models.

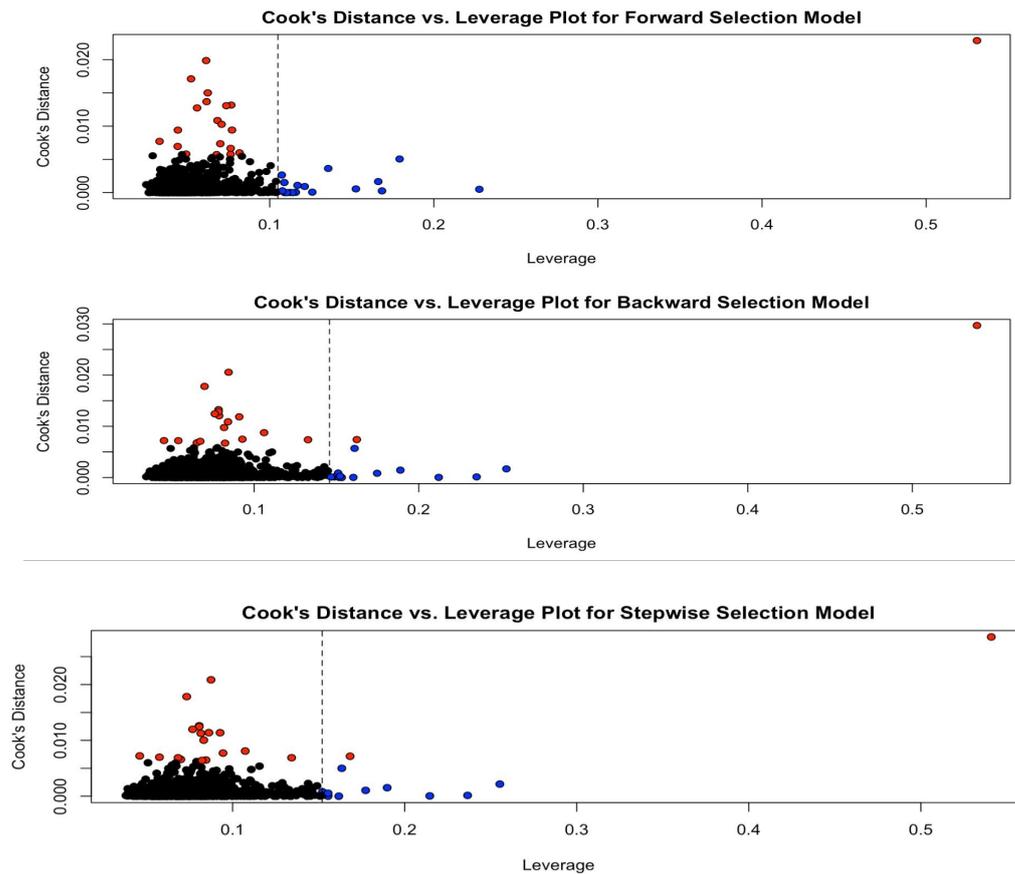

**Figure 5. Leverage vs. Cook's Distance Plot for all Three Models**

As we observed in the exploratory data analysis, observation 985 is an outlier with extremely high leverage and influence compared to other observations. We can see that among the 15 observations with the highest influence (labelled in red), there are 1 observation in the backward and stepwise model that has both high leverage and high influence compared to 0 in the forward selection model (excluding observation 985).

We also noticed that although there are less high leverage observations in the backward and stepwise model (labelled in blue), there are more observations with high influences in the backward and stepwise model, since the observations are more spread out along the x-axis compared to the forward selection model.

According to the scatterplot above, observation 985 has more influence on the backward and stepwise model (the Cook's distance for observation 985 is about 0.030 for backward model and 0.029 for the stepwise model, compared to around 0.024 for forward selection model). Lastly, we will perform cross-validation that runs for 8000 replications using 80% of the observations for the training set for good accuracy.

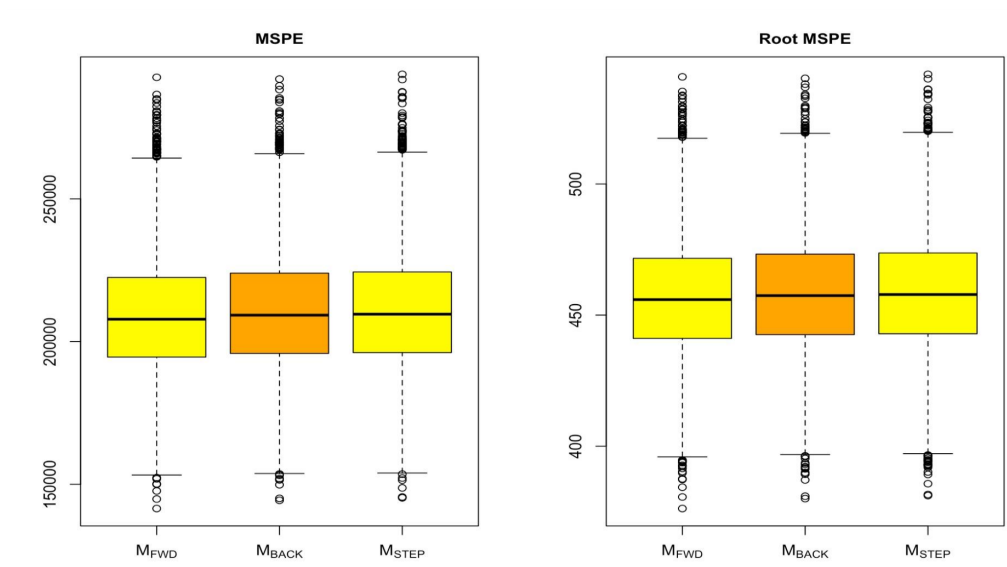

**Figure 6. Cross-Validation with 8000 Repetitions and Training Size of 80% Sample Size**

From Figure 6, we obtain box plots followed by their interquartile range information. This data tells how narrow are the ranges between the first and the third quantile in our accuracy values (with IQR). From the observations over the above data, we see that the forward selection model has the narrowest IQR (27892), which means that 50% of the values are found in a narrower range compared to the other two models. Meanwhile, considering all data, it is still the forward model showing all data laying in the narrowest full range (between minimum and maximum). Therefore, we say that the values of the forward selection model are the most condensed, and thus it is the best fitted model. The detailed five-number summary for both boxplots could be found in Appendix B.

Lastly, we use the dataset without the outlier observation 985 to fit all of the forward, backward and stepwise models to see if our conclusions are substantively different.

```
                           Forward       Backward       Stepwise
sum-of-squared PRESS  262641303.3298 262081477.7022 262081477.7022
AIC                        19216.7657     19215.4908     19215.4908
Adjusted R^2                   0.2448         0.2575         0.2575
sum-of-squared DFFITS         76.5534        99.2822        99.2822
Number of Predictors          72.0000        94.0000        94.0000
```

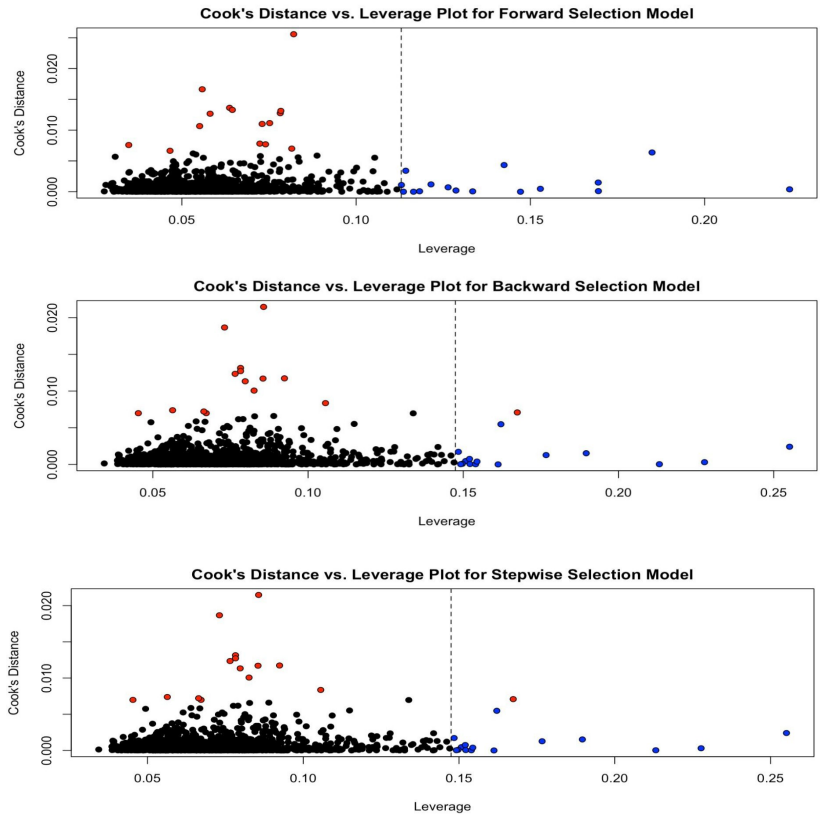

**Figure 7. Leverage vs. Cook's Distance for All Three Refitted Models Using predict_data[-985, ]**

We can see that the results slightly change since our Backward selection model is preferred with adjusted R squared, AIC and sum-of-squared PRESS. We also observed that there are significantly less high-leverage observations in the backward elimination model, and the observation that has the highest Cook's distance are not as far from the others as they are in the forward selection model. However, the leverage of the observation within the backward and stepwise model are still generally larger than the observations in the forward selection model, since the observation points are more horizontally spreaded out along the x-axis and the cutoff of high-leverage observations for backward and stepwise model is near 0.15, while that of forward selection model is around 0.11.

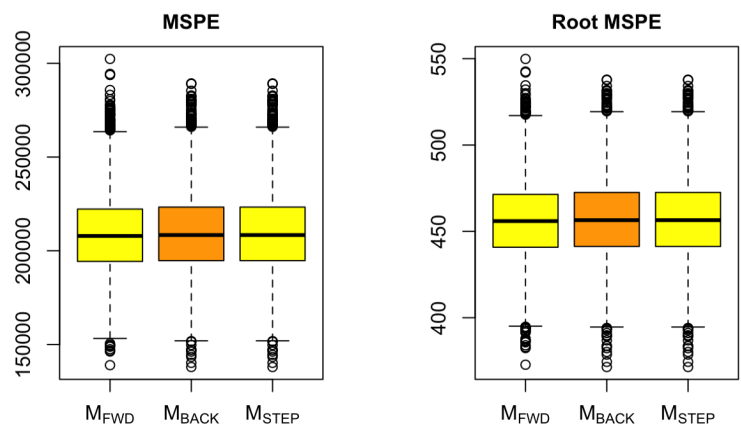

**Figure 8. Cross-Validation with 8000 Repetitions and Training Size of 80% Sample Size**

Also, by observing the boxplot after performing the cross validation that runs for 8000 replications, the forward selection model shows some very large MSPEs (see the above boxplot in Figure 8) and it is more likely to overfit the data than the backward and stepwise model. However, we can see that the IQR for the forward selection is the narrowest (27926.5). Similar to the analysis towards Figure 6, this means that the forward selection still demonstrates the most condensed data. Thus, we also conclude that forward selection is the best model here. The detailed five-number summary for both boxplots could be found in Appendix B.

## Results

According to the above analysis, we pick the forward selection model as our final selected model. The reasons is as follows:
- First, although the backward and stepwise model has slightly higher adjusted R^2 value, the forward selection model has the advantage of being more easily interpretable by a large margin, since the number of predictors in the forward selection model is 67 in contrast to 93 and 97 in the backward selection model.
- Second, the observations in the forward selection model has lower leverage in general compared to the backward and stepwise selection model as shown by the Leverage vs. Cook's Distance scatterplot. As a result, the forward selection model has significantly lower sum-of-squared DFFITS values compared to the other two models.
- Also, from our observations towards both box plots (with and without outliers), the forward selection model has the narrowest IQR, having all data being more condensed together.

We proceed to examine whether our selected model satisfies the four assumptions of multiple linear regression.

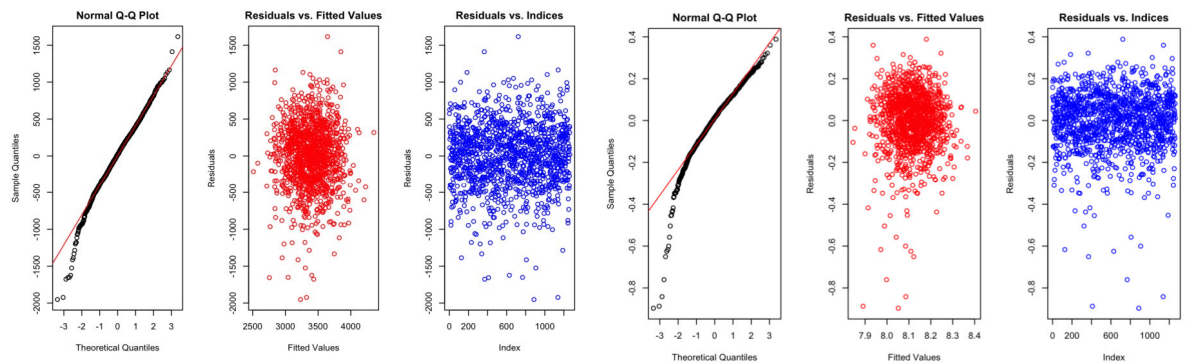

**Figure 9. Residual Plots Summary**
**Without Log Transformation (left) and Without Log Transformation (Right) on the Response Variable**

Since both the scatterplot of Residuals vs Fitted Values shows a violation of the homoscedasticity assumption, we consider the log transformation on the response variable e3_bw. After applying the log transformation on the response variable e3_bw, we get the following result on the right side of Figure 9. After the log transformation on the response variable, we can see that from the QQ plot, we can see that most residuals lie on the theoretical normal qqline except for some observations near the two tails. Also, we can see that from the histogram of studentized residuals, the residuals look like the density for N(0, 1). This suggests that our selected model **satisfies the normality assumption of MLR.**

For the scatterplot between residuals and indices, the observations are evenly spread out around 0, which suggests that the residuals of our selected model are independent. Considering how the data was collected, since the 1301 observations of mother-child pairs was being collected from hospitals, clinics and parental care centres of six different countries, it is likely that observations will cluster within hospitals since mother-child pairs from the same city would have similar chemical exposures. However, since each mother-child pair have different habits and living environment, it is very likely that these observations do not depend on each other. Hence, this suggests our selected model **satisfies the independence assumption of MLR.**

Observing the plot Residuals vs. Fitted Values, we can see that the majority of the points lie on the line y = 0. Although there are some outliers in the middle of the graph, since we have a significant amount of values that form a horizontal band around the line y = 0, it suggests that our selected model **satisfies the equal variance assumption of MLR.**

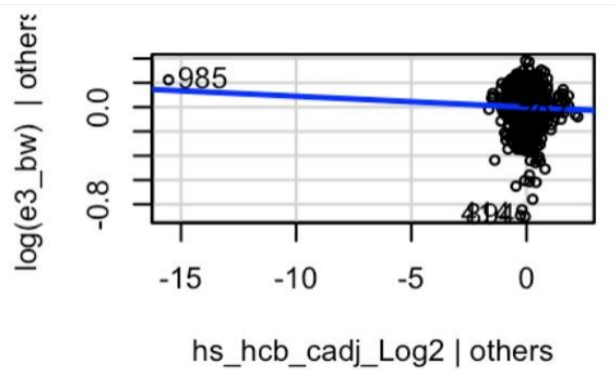

Figure 11. Added-Variate Plot

By observing the Added-Variable plots for our selected data, we can see that the covariates we decided to keep in our final model do meet the linearity assumptions. Although the graphs for some covariates may seem to be questionable in terms of being linear to the response, however if we take a look closely, the scale of the graph can make the visual of the outcomes misleading. For example, the plot for the covariate "hs_hcb_cadj_Log2" in the figure 11 above, the points may seem to form a vertical relationship where the line of best fit is almost horizontal, this is because the scale of the y-axis is in increments of 0.2 units, but the x-axis is in increments of 5 units due to an extreme outlier. Since the x-axis is in units larger than the y-axis, this makes the horizontal spread of the points in the graph less obvious than the vertical spread of the points. As a result, it suggests that our model **satisfies the linear assumption of MLR.**

# Discussions

In summary, the results of our calculations show that by eliminating the covariates with a significant amount of N/A entries, then omitting rows with one or more N/A entries, and finally applying a forward selection on the reduced dataset and applying a logarithmic transformation on the respons variable, we have reached a best fitting model for birthweights. Under testing, the results show that the covariates meet the assumptions of linearity, normality, homoscedasticity, and independence.

Based on the results of our findings, building a model with forward selection and adding a log transformation to the response has created the most fitting predictive model for birthweight. This model has met all the assumptions of MLR. Furthermore, we have cross validated that using forward selection is more suitable than other selection methods such as backward and stepwise selection.

The limitation of our observation is from our way of handling missing values. We drop some whole columns for their large numbers of missing values in the certain category. Dropping these categories helps to focus on the existing data if the missing values are not important for our prediction. However, if these given categories are important, then dropping relevant categories could result in lack of observations and factors for a strong model. The harm of removing these categories manually is that such loss of information is not at random. Specific sections of observations could be accidentally removed and ignored by these decisions. Some possible ways to improve our limitation could be concluding missing values from mathematical ways or logical ways. For example, we could consider the mean, median, or mode of the existing data to deduce the missing values, or simply guess by considering what makes more sense. These assumptions and deductions allow the potential effect of the incomplete variables to the model. Taking these variables into consideration could improve our analysis to be more well-rounded in terms of considerations towards all factors.

# Appendix

## Appendix A (Final Selected Model)
> summary(Mfwd)

Call:
lm(formula = e3_bw ~ h_cohort + hs_c_height_None + e3_sex_None +
    hs_mo_c_Log2 + h_Benzene_Log + h_folic_t1_None + hs_pcb118_cadj_Log2 +
    hs_cs_m_Log2 + hs_dmtp_madj_Log2 + hs_mibp_madj_Log2 + e3_asmokcigd_p_None +
    hs_mbzp_cadj_Log2 + h_builtdens300_preg_Sqrt + hs_no2_dy_hs_h_Log +
    hs_mepa_madj_Log2 + hs_dep_madj_Log2 + hs_pfoa_c_Log2 + hs_c_weight_None +
    hs_pfoa_m_Log2 + hs_dde_madj_Log2 + hs_total_sweets_Ter +
    hs_pbde153_madj_Log2 + h_NO2_Log + hs_pet_dog_r2_None + hs_mnbp_madj_Log2 +
    hs_fastfood_Ter + h_accesspoints300_preg_Log + hs_accesspoints300_h_Log +
    hs_hg_c_Log2 + hs_org_food_Ter + h_temperature_preg_None +
    h_Absorbance_Log + hs_blueyn300_h_None + hs_pm25abs_wk_hs_h_Log +
    hs_hg_m_Log2 + hs_bupa_cadj_Log2 + hs_accesslines300_s_dic0 +
    hs_accesslines300_h_dic0 + hs_popdens_s_Sqrt + hs_greenyn300_s_None +
    hs_hcb_madj_Log2 + hs_mep_madj_Log2 + h_walkability_mean_preg_None +
    h_connind300_preg_Sqrt + hs_connind300_s_Log + hs_dde_cadj_Log2 +
    hs_hcb_cadj_Log2 + hs_tl_cdich_None + hs_pcb153_madj_Log2 +
    hs_pcb170_madj_Log2 + h_PM_Log + hs_dif_hours_total_None +
    hs_cd_m_Log2 + hs_pb_m_Log2 + h_lden_cat_preg_None +
h_landuseshan300_preg_None +
    hs_landuseshan300_s_None + hs_trafload_h_pow1over3 + hs_fdensity300_h_Log,
    data = predict_data)

Residuals:
     Min       1Q   Median       3Q      Max
-1966.71  -265.85    18.02   291.26  1598.56

Coefficients:
                      Estimate  Std. Error  t value  Pr(>|t|)
(Intercept)           1929.6357   516.5647    3.736  0.000196 ***
h_cohort2             -305.1608   110.1154   -2.771  0.005669 **
h_cohort3             -318.0075   109.6092   -2.901  0.003784 **
h_cohort4              104.8506    98.0250    1.070  0.284999
h_cohort5              -24.5711    77.8106   -0.316  0.752223
h_cohort6             -161.0267   114.3618   -1.408  0.159374
hs_c_height_None      1203.6412   277.2396    4.342  1.53e-05 ***
e3_sex_Nonemale        105.5327    26.0630    4.049  5.47e-05 ***
hs_mo_c_Log2           -45.8800    14.4452   -3.176  0.001530 **

```
h_Benzene_Log            -138.5237   38.4712 -3.601 0.000330 ***
h_folic_t1_None1          -82.2791   30.3376 -2.712 0.006780 **
hs_pcb118_cadj_Log2        58.7589   21.5838  2.722 0.006575 **
hs_cs_m_Log2              -91.1216   29.0337 -3.138 0.001739 **
hs_dmtp_madj_Log2          14.2454    4.3950  3.241 0.001222 **
hs_mibp_madj_Log2          38.3324   13.7234  2.793 0.005301 **
e3_asmokcigd_p_None       -29.4843    8.0202 -3.676 0.000247 ***
hs_mbzp_cadj_Log2          27.0323   11.2735  2.398 0.016642 *
h_builtdens300_preg_Sqrt   -0.2523    0.1299 -1.942 0.052374 .
hs_no2_dy_hs_h_Log         52.9976   27.0989  1.956 0.050730 .
hs_mepa_madj_Log2          -9.3437    5.4524 -1.714 0.086847 .
hs_dep_madj_Log2          -24.2391    8.4062 -2.883 0.004003 **
hs_pfoa_c_Log2             53.9611   25.5053  2.116 0.034576 *
hs_c_weight_None            9.0469    2.8790  3.142 0.001716 **
hs_pfoa_m_Log2            -31.5235   14.5248 -2.170 0.030176 *
hs_dde_madj_Log2          -21.1021    8.1018 -2.605 0.009311 **
hs_total_sweets_Ter(4.1,8.5] 23.2941 31.9309  0.730 0.465828
hs_total_sweets_Ter(8.5,Inf] -63.0624 35.1250 -1.795 0.072844 .
hs_pbde153_madj_Log2       12.2231    4.8340  2.529 0.011580 *
h_NO2_Log                  39.8297   17.0994  2.329 0.020007 *
hs_pet_dog_r2_None1        76.7645   37.0057  2.074 0.038253 *
hs_mnbp_madj_Log2         -29.4947   12.1116 -2.435 0.015025 *
hs_fastfood_Ter(0.132,0.5] -120.4555 43.3275 -2.780 0.005518 **
hs_fastfood_Ter(0.5,Inf]  -111.5417  44.7350 -2.493 0.012786 *
h_accesspoints300_preg_Log 104.8261  31.3240  3.347 0.000844 ***
hs_accesspoints300_h_Log  -60.3310   22.4857 -2.683 0.007394 **
hs_hg_c_Log2              -20.1983    9.1657 -2.204 0.027735 *
hs_org_food_Ter(0.132,1]   23.6631   35.3064  0.670 0.502843
hs_org_food_Ter(1,Inf]     82.9697   35.4028  2.344 0.019260 *
h_temperature_preg_None     9.5431    5.5415  1.722 0.085305 .
h_Absorbance_Log          101.8770   37.2628  2.734 0.006348 **
hs_blueyn300_h_None1       79.7307   44.6638  1.785 0.074491 .
hs_pm25abs_wk_hs_h_Log    -81.9019   39.0480 -2.097 0.036159 *
hs_hg_m_Log2               22.8289   11.2497  2.029 0.042647 *
hs_bupa_cadj_Log2         -12.8339    6.4829 -1.980 0.047969 *
hs_accesslines300_s_dic01 -229.1696  71.2012 -3.219 0.001322 **
hs_accesslines300_h_dic01  172.5072  82.4585  2.092 0.036642 *
hs_popdens_s_Sqrt           1.0610    0.4359  2.434 0.015088 *
hs_greenyn300_s_None1      78.2037   33.7728  2.316 0.020748 *
hs_hcb_madj_Log2           28.5745   13.6363  2.095 0.036337 *
hs_mep_madj_Log2           13.2552    7.4925  1.769 0.077125 .
h_walkability_mean_preg_None -771.1232 262.7191 -2.935 0.003397 **
h_connind300_preg_Sqrt     10.1116    4.2886  2.358 0.018543 *
hs_connind300_s_Log       -49.9940   22.0852 -2.264 0.023770 *
```

```
hs_dde_cadj_Log2              25.4545   12.1498   2.095 0.036374 *
hs_hcb_cadj_Log2             -36.3590   20.1118  -1.808 0.070879 .
hs_tl_cdich_NoneUndetected   -78.6128   48.4199  -1.624 0.104730
hs_pcb153_madj_Log2          -65.3800   24.1634  -2.706 0.006911 **
hs_pcb170_madj_Log2           38.1222   17.8000   2.142 0.032417 *
h_PM_Log                     -56.0102   31.4307  -1.782 0.074997 .
hs_dif_hours_total_None      -33.0016   20.6742  -1.596 0.110689
hs_cd_m_Log2                  27.5982   15.4293   1.789 0.073915 .
hs_pb_m_Log2                 -31.8619   19.0069  -1.676 0.093931 .
h_lden_cat_preg_None           3.3006    2.0324   1.624 0.104628
h_landuseshan300_preg_None   234.2263  119.0658   1.967 0.049388 *
hs_landuseshan300_s_None    -227.7371  129.6853  -1.756 0.079329 .
hs_trafload_h_pow1over3       -0.4363    0.2896  -1.506 0.132216
hs_fdensity300_h_Log          25.9275   17.8415   1.453 0.146423
---
Signif. codes:  0 '***' 0.001 '**' 0.01 '*' 0.05 '.' 0.1 ' ' 1

Residual standard error: 442.2 on 1209 degrees of freedom
Multiple R-squared:  0.2827,  Adjusted R-squared:  0.2435
F-statistic: 7.22 on 66 and 1209 DF,  p-value: < 2.2e-16
```

**Appendix B (Resulted MSPE from the Cross-Validation for all three models by using the dataset predict_data and predict_data[-985, ])**

Summary of MSPE and squared MSPE for all Three Models Fitted Using predict_data:

```
          MSPE_Mfwd MSPE_Mback MSPE_Mstep
Min.         141521     144386     145269
1st Qu.      194558     195846     196141
Median       207840     209252     209611
Mean         209022     210396     210650
3rd Qu.      222450     223980     224377
Max.         292560     291954     293621
          MSPE_sqr_Mfwd MSPE_sqr_Mback MSPE_sqr_Mstep
Min.              376.2          380.0          381.1
1st Qu.           441.1          442.5          442.9
Median            455.9          457.4          457.8
Mean              456.6          458.1          458.4
3rd Qu.           471.6          473.3          473.7
Max.              540.9          540.3          541.9
```

Summary of MSPE and squared MSPE for all Three Models Fitted Using predict_data[-985, ]:

```
          MSPE_Mfwd MSPE_Mback MSPE_Mstep
Min.        138980     137977     137977
1st Qu.     194322     194743     194743
Median      207904     208390     208390
Mean        208760     209334     209334
3rd Qu.     222248     223289     223289
Max.        302344     289310     289310
          MSPE_sqr_Mfwd MSPE_sqr_Mback MSPE_sqr_Mstep
Min.              372.8          371.5          371.5
1st Qu.           440.8          441.3          441.3
Median            456.0          456.5          456.5
Mean              456.3          457.0          457.0
3rd Qu.           471.4          472.5          472.5
Max.              549.9          537.9          537.9
```

**Appendix C (R code)**

```r
###############################################################################
######################## DATA PREPARATION #####################

# Import libraries
library(ggplot2)
library(corrplot)
library(dplyr)
library(car)

knitr::opts_chunk$set(echo = TRUE)
load("exposome_NA.RData")
phenotypes <- phenotypeNA
covariates <- covariatesNA
exposome <- exposomeNA[order(exposomeNA$ID),]

# The maximum number of observations that can be eliminated
NA_remove_ratio <- 0.01
covariate_rows <- nrow(covariates)   # number of rows in the covariate table
exposome_rows <- nrow(exposome) # number of rows in the exposome table

reduced_covariates <- data.frame(covariates)
reduced_exposome <- data.frame(exposome)

for (i in 2:ncol(covariates)) {
  if (sum(is.na(covariates[ , i])) >= covariate_rows * NA_remove_ratio) {
```

```r
    var_del <- names(covariates)[i]
    reduced_covariates <-reduced_covariates[,-(which(names(reduced_covariates) %in% var_del))]
  }
}

for (i in 2:ncol(exposome)) {
  if (sum(is.na(exposome[ , i])) >= exposome_rows * NA_remove_ratio)
  { var_del <- names(exposome)[i]
    reduced_exposome <-reduced_exposome[,-(which(names(reduced_exposome) %in% var_del))]
  }
}

predict_data <- merge(reduced_covariates, reduced_exposome, by = "ID")
predict_data$e3_bw <- phenotypes$e3_bw
predict_data <- na.omit(predict_data)

# transform some of the variables from numeric to factor for
predict_data$h_accesslines300_preg_dic0 <- factor(predict_data$h_accesslines300_preg_dic0)
predict_data$hs_accesslines300_h_dic0 <- factor(predict_data$hs_accesslines300_h_dic0)
  predict_data$hs_accesslines300_s_dic0 <- factor(predict_data$hs_accesslines300_s_dic0)

    ###############################################################################
  ##################### EXPLORATORY DATA ANALYSIS ############################

numeric_rows <- nrow(predict_data)  ## total number of observations
numeric_data <- data.frame(predict_data)
categorical_data <- data.frame(predict_data)

#Extract all numeric data into one dataset, and all categorical data in another
for (i in 2:ncol(predict_data)) {
  if (is.numeric(predict_data[,i]) == FALSE)
   { var_del <- names(predict_data)[i]
     numeric_data <-numeric_data[,-(which(names(numeric_data) %in% var_del))] } else {
     var_del <- names(predict_data)[i]
     categorical_data <- categorical_data[,-(which(names(categorical_data) %in% var_del))] }
}

# produce scatterplot matrices
pairs(numeric_data[,c(11:20)], main = "Figure 1:
```

Paired Data between the 11th and 20th Numeric Variables", cex = 0.2)

```r
## compute leverage
M <- lm(e3_bw~., data = predict_data)
lev <- hatvalues(M) ## leverage (h_i)
hbar <- mean(lev) ## \bar{h}
c(sum(lev),ncol(model.matrix(M))) ## check trace is same as rank of
lev_high <- which(lev > 2 * hbar) ## x values for labelling points > 2hbar
#ids <- which(lev>2*hbar) ## x values for labelling points >2hbar

## compute Cook's Distance
cook <- cooks.distance(M)
cook_high <- (cook >= quantile(cook,
                    probs = ((nrow(predict_data) - 10) / nrow(predict_data))))
cook_high_ids <- which(cook >= quantile(cook,
                    probs = ((nrow(predict_data) - 10) / nrow(predict_data))))

## plot leverage for the forward selection model
colors <- rep("black", len = nrow(predict_data))
colors[lev_high] <- "blue"
colors[cook_high] <- "red"
plot(lev, cook, xlab = "Leverage", ylab = "Cook's Distance",
     main = "Cook's Distance vs. Leverage Plot on the Full Model", pch = 21, bg = colors)
abline(v = 2 * hbar, lty = 2, col = "purple")

## label points > 2hbar
text(x=lev[lev_high],y=cook[lev_high], labels=lev_high, cex= 0.6, pos=2)
## label points > 2hbar
text(x=lev[cook_high_ids],y=cook[cook_high_ids], labels=cook_high_ids, cex= 0.6, pos=2)
# add legend
legend("topright", legend = c("Top 10 points that
have Cook's Distance > 0.5", "leverage > 2 * mean(Leverage)",
                   "2 * mean(Leverage)"),
       col = c("red", "blue", "purple"), cex = 0.7, pch = c(19, 19, 19, NA), lty = c(NA, NA, NA, 2), bty = "n")
```

METHODS

```
###############################################################################
############### ELIMINATE COVARIATES WITH HIGH VIFs ##################
```

```r
reduce_data <- function(vstar){
  mod <- lm(e3_bw ~ ., data = numeric_data)
  df_temp <- numeric_data
  while(1){
    mod <- lm(e3_bw ~ ., data = df_temp) # fit the model on the reduced dataframe
    if(max(car::vif(mod)) > vstar){
      var_del <- names(which.max(car::vif(mod)))
      print(var_del)
      # extract the name of the variable that has the highest vif
      ## then remove this variable from the dataframe
      df_temp <- df_temp[,-(which(names(df_temp) %in% var_del))]
    } else { # if the maximum vif among the current covariates is small (i.e. <= vstar)
      break
    }
  }
  return(df_temp)
}

# get the set of all numeric data in the predict_data dataset
numeric_reduced <- reduce_data(10)
predict_data <- merge(numeric_reduced, categorical_data, by = "ID")

# plot histogram
hist(vif(lm(e3_bw ~ ., data = numeric_data)), xlab = "VIF",
     main = "VIF Distribution before Elimination", col = "cyan")

hist(vif(lm(e3_bw ~ ., data = numeric_reduced)), xlab = "VIF",
     main = "VIF Distribution after Elimination", col = "cyan")

####################################################################################
##################### AUTOMATED MODEL SELECTION #####################

predict_data_remove985 <- predict_data[-985, ]

M0 <- lm(e3_bw ~ 1, data = predict_data)  ## minimal model

Mfull <- lm(e3_bw ~ ., data = predict_data) ## full model

df.penalty <- 2 ## this is the k penalty

# forward selection
```

```r
system.time({
Mfwd <- step(object = M0, # base model
        scope = list(lower = M0, upper = Mfull), # smallest and largest model
        trace = 0, # trace prints out information
        direction = "forward", k =
        df.penalty) })

# backward selection
system.time({
  Mback <- step(object = Mfull, # base model
        scope = list(lower = M0, upper = Mfull),
        direction = "backward", trace = 0)
})

# stepwise (both directions) selection
Mstart <- lm(e3_bw ~ ., data = predict_data) # starting point model: main effects only
system.time({
  Mstep <- step(object = Mstart,
        scope = list(lower = M0, upper = Mfull),
        direction = "both", trace = 0, k = df.penalty)
})

# combining results
models <- c(Mfwd$call, Mback$call, Mstep$call)
names(models) <- c("FWD", "BACK", "STEP")
models

####################### Refit all three models after removing observation 985 #######################

predict_data <- predict_data[-985, ]

M0 <- lm(e3_bw ~ 1, data = predict_data)  ## minimal model

Mfull <- lm(e3_bw ~ ., data = predict_data) ## full model

df.penalty <- 2 ## this is the k penalty

# forward selection
system.time({
Mfwd_remove985 <- step(object = M0, # base model
        scope = list(lower = M0, upper = Mfull), # smallest and largest model
        trace = 0, # trace prints out information
        direction = "forward", k =
        df.penalty) })
```

```r
# backward selection
system.time({
  Mback_remove985 <- step(object = Mfull, # base model
          scope = list(lower = M0, upper = Mfull),
          direction = "backward", trace = 0)
})

# stepwise (both directions) selection
Mstart <- lm(e3_bw ~ ., data = predict_data) # starting point model: main effects only
system.time({
  Mstep_remove985 <- step(object = Mstart,
          scope = list(lower = M0, upper = Mfull),
          direction = "both", trace = 0, k = df.penalty)
})

    ###########################################################################
    ##################### MODEL DATA COMARISON TABLE #####################

# Disable scientific notation, set to display 4 decimal places
options(scipen = 999, digits=4)

## Define
m1 <- Mfwd
m2 <- Mback
m3 <- Mstep
m1.h <- hatvalues(m1)
m2.h <- hatvalues(m2)
m3.h <- hatvalues(m3)
m1.PRESS <- resid(m1)/(1-m1.h)
m2.PRESS <- resid(m2)/(1-m2.h)
m3.PRESS <- resid(m3)/(1-m3.h)
m1.PRESS.sq <- sum(m1.PRESS^2)  ## sum of square of PRESS
m2.PRESS.sq <- sum(m2.PRESS^2)
m3.PRESS.sq <- sum(m3.PRESS^2)

# AIC
m1.AIC <- AIC(m1)
m2.AIC <- AIC(m2)
m3.AIC <- AIC(m3)

# Adjusted R^2
```

```r
m1.R2 <- summary(m1)$adj.r.squared
m2.R2 <- summary(m2)$adj.r.squared
m3.R2 <- summary(m3)$adj.r.squared

# sum of square DFFITS
m1.dffits.sq <- sum(dffits(m1)^2)
m2.dffits.sq <- sum(dffits(m2)^2)
m3.dffits.sq <- sum(dffits(m3)^2)

# Number of Betas
m1.rank <- m1$rank
m2.rank <- m2$rank
m3.rank <- m3$rank

# Construct a comparison table
## combining them
PRESS.list <- c(m1.PRESS.sq, m2.PRESS.sq, m3.PRESS.sq)
AIC.list <- c(m1.AIC, m2.AIC, m3.AIC)
R2.list <- c(m1.R2, m2.R2, m3.R2)
DFFITS.list <- c(m1.dffits.sq, m2.dffits.sq, m3.dffits.sq)
rank.list <- c(m1.rank, m2.rank, m3.rank)
diagnost.matrix <- rbind(PRESS.list, AIC.list, R2.list, DFFITS.list, rank.list)
row.names(diagnost.matrix) <- c("sum-of-squared PRESS", "AIC", "Adjusted R^2",
                    "sum-of-squared DFFITS", "Number of Predictors")
colnames(diagnost.matrix) <- c("Forward", "Backward", "Stepwise")

# Output the table
diagnost.matrix

###############################################################################
We will perform the analyzation again with refitted models after deleting the outlier observation 985,
by setting m1 <- Mfwd_remove985, m2 <- Mback_remove985, m3 <- Mstep_remove985

###############################################################################

    ###############################################################################
        ############# COOK'S DISTANCE vs. LEVERAGE SCATTERPLOT #############

## compute leverage
M1 <- Mfwd
M2 <- Mback
```

```r
M3 <- Mstep
## leverage (h_i)
lev1 <- hatvalues(M1)
lev2 <- hatvalues(M2)
lev3 <- hatvalues(M3)
## compute h bar
hbar1 <- mean(lev1)
hbar2 <- mean(lev2)
hbar3 <- mean(lev3)
## check trace is same as rank
c(sum(lev1),ncol(model.matrix(M1)))
c(sum(lev2),ncol(model.matrix(M2)))
c(sum(lev3),ncol(model.matrix(M3)))
## x values for labelling points > 2hbar
lev1_high <- which(lev1 > 2 * hbar1)
lev2_high <- which(lev2 > 2 * hbar2)
lev3_high <- which(lev3 > 2 * hbar3)

# compute Cook's distance
cook1 <- cooks.distance(M1)
cook2 <- cooks.distance(M2)
cook3 <- cooks.distance(M3)
cook1_high <- (cook1 >= quantile(cook1,
                    probs = ((nrow(predict_data) - 15) / nrow(predict_data))))
cook2_high <- (cook2 >= quantile(cook2,
                    probs = ((nrow(predict_data) - 15) / nrow(predict_data))))
cook3_high <- (cook3 >= quantile(cook3,
                    probs = ((nrow(predict_data) - 15) / nrow(predict_data))))

par(mfrow = c(2, 1))
## plot leverage for the forward selection model
colors1 <- rep("black", len = nrow(predict_data))
colors1[lev1_high] <- "blue"
colors1[cook1_high] <- "red"
plot(lev1, cook1, xlab = "Leverage", ylab = "Cook's Distance",
    main = "Cook's Distance vs. Leverage Plot for Forward Selection Model", pch = 21, bg = colors1)
abline(v = 2 * hbar1, lty = 2)

## plot leverage for the backward selection model
colors2 <- rep("black", len = nrow(predict_data))
colors2[lev2_high] <- "blue"
colors2[cook2_high] <- "red"
```

```r
plot(lev2, cook2, xlab = "Leverage", ylab = "Cook's Distance",
     main = "Cook's Distance vs. Leverage Plot for Backward Selection Model", pch = 21, bg = colors2)
abline(v = 2 * hbar2, lty = 2)

par(mfrow = c(2, 1))
## plot leverage for the backward selection model
colors3 <- rep("black", len = nrow(predict_data))
colors3[lev3_high] <- "blue"
colors3[cook3_high] <- "red"
plot(lev3, cook3, xlab = "Leverage", ylab = "Cook's Distance",
     main = "Cook's Distance vs. Leverage Plot for Stepwise Selection Model", pch = 21, bg = colors3)
abline(v = 2 * hbar3, lty = 2)

################################################################################
We will perform the analyzation again with refitted models after deleting the outlier observation 985,
by setting m1 <- Mfwd_remove985, m2 <- Mback_remove985, m3 <- Mstep_remove985

################################################################################

     ################################################################################
     ##################### CROSS VALIDATION AND BOXPLOT #####################

set.seed(20883271)

# compare Mfull, Mback and Mstep (cross validation)
M1 <- Mfwd
M2 <- Mback
M3 <- Mstep
numeric_rows <- nrow(predict_data)   # total number of observation
num_train <- round(0.8 * numeric_rows) # use 80% of the observation as the training set
num_test <- numeric_rows-num_train   # size of testing set
nreps <- 8000                        # number of replications
# sum of error for cross validation replication
sse1 <- rep(NA, nreps)
sse2 <- rep(NA, nreps)
sse3 <- rep(NA, nreps)

mspe1 <- rep(NA, nreps)
mspe2 <- rep(NA, nreps)
mspe3 <- rep(NA, nreps)
```

```r
# measure time
system.time({
  for( ii in 1:nreps ) {
    if ( ii%%100 == 0 ) message ( "ii = " , ii)
    
    # randomly select training observations
    train.ind <- sample(numeric_rows, num_train)
    M1.cv <- update(M1, subset = train.ind)
    M2.cv <- update(M2, subset = train.ind)
    M3.cv <- update(M3, subset = train.ind)
    
    # test models
    M1.res <- predict_data$e3_bw[-train.ind] -
      predict(M1.cv, newdata = predict_data[-train.ind,])
    M2.res <- predict_data$e3_bw[-train.ind]-
      predict(M2.cv, newdata = predict_data[-train.ind,])
    M3.res <- predict_data$e3_bw[-train.ind]-
      predict(M3.cv, newdata = predict_data[-train.ind,])
    
# Calculate MSPEs
    mspe1[ii] <- mean(M1.res^2)
    mspe2[ii] <- mean(M2.res^2)
    mspe3[ii]                  <-
    mean(M3.res^2) }})

# compare
Mnames <- expression(M[FWD], M[BACK], M[STEP])
par(mfrow = c(1,2))
cex <- 1
boxplot(x = list(mspe1, mspe2, mspe3), names = Mnames,
     main = "MSPE",
     #ylab = expression(sqrt(bar(SSE)[CV])),
     #ylab = expression(MSPE),
     col = c("yellow", "orange"),
     cex = cex, cex.lab = cex, cex.axis = cex, cex.main = cex)
boxplot(x = list(sqrt(mspe1), sqrt(mspe2), sqrt(mspe3)), names = Mnames,
     main = "Root MSPE",
     #ylab = expression(sqrt(MSPE)),
     # ylab = expression(SSE[CV]),
     col = c("yellow", "orange"),
     cex = cex, cex.lab = cex, cex.axis = cex, cex.main = cex)
```

```r
# Output the MSPE summary table for all three models
cbind(MSPE_Mfwd=summary(mspe1),
MSPE_Mback=summary(mspe2),
MSPE_Mstep=summary(mspe3))

# Output the squared MSPE summary table for all three models
cbind(MSPE_sqr_Mfwd=summary(sqrt(mspe1)),
    MSPE_sqr_Mback=summary(sqrt(mspe2)),
MSPE_sqr_Mstep=summary(sqrt(mspe3)))

###############################################################################
We will perform the analyzation again with refitted models after deleting the outlier observation 985,
by setting m1 <- Mfwd_remove985, m2 <- Mback_remove985, m3 <- Mstep_remove985

###############################################################################

    ##################################################################################
   ################### MLR ASSUMPTION CHECK FOR THE SELECTED MODEL############

summary(Mfwd)$call
# call summary(Mfwd)$call to get the formula first, then transform the response variate using
log transformation
Mfwd_transformed <- lm(formula = log(e3_bw) ~ h_cohort + hs_c_height_None +
e3_sex_None +
    hs_mo_c_Log2 + h_Benzene_Log + h_folic_t1_None + hs_pcb118_cadj_Log2 +
    hs_cs_m_Log2 + hs_dmtp_madj_Log2 + hs_mibp_madj_Log2 + e3_asmokcigd_p_None
+
    hs_mbzp_cadj_Log2 + h_builtdens300_preg_Sqrt + hs_no2_dy_hs_h_Log +
    hs_mepa_madj_Log2 + hs_dep_madj_Log2 + hs_pfoa_c_Log2 + hs_c_weight_None +
    hs_pfoa_m_Log2 + hs_dde_madj_Log2 + hs_total_sweets_Ter +
    hs_pbde153_madj_Log2 + h_NO2_Log + hs_pet_dog_r2_None + hs_mnbp_madj_Log2 +
    hs_fastfood_Ter + h_accesspoints300_preg_Log + hs_accesspoints300_h_Log +
    hs_hg_c_Log2 + hs_org_food_Ter + h_temperature_preg_None +
    h_Absorbance_Log + hs_blueyn300_h_None + hs_pm25abs_wk_hs_h_Log +
    hs_hg_m_Log2 + hs_bupa_cadj_Log2 + hs_accesslines300_s_dic0 +
    hs_accesslines300_h_dic0 + hs_popdens_s_Sqrt + hs_greenyn300_s_None +
    hs_hcb_madj_Log2 + hs_mep_madj_Log2 + h_walkability_mean_preg_None +
    h_connind300_preg_Sqrt + hs_connind300_s_Log + hs_dde_cadj_Log2 +
    hs_hcb_cadj_Log2 + hs_tl_cdich_None + hs_pcb153_madj_Log2 +
```

```r
    hs_pcb170_madj_Log2 + h_PM_Log + hs_dif_hours_total_None +
    hs_cd_m_Log2 + hs_pb_m_Log2 + h_lden_cat_preg_None + h_landuseshan300_preg_None +
    hs_landuseshan300_s_None + hs_trafload_h_pow1over3 + hs_fdensity300_h_Log,
    data = predict_data)

M_selected <- Mfwd_transformed

# plot normal qqplot for residuals
par(mfrow = c(1, 3))
qqnorm(resid(M_selected))
qqline(M_selected$residuals, col = "red")

# Check Normality Assumptions
plot(M_selected$residuals~M_selected$fitted.values,
     main = "Residuals vs. Fitted Values", ylab = "Residuals", xlab = "Fitted Values",
     col = "red")
plot(M_selected$residuals, main = "Residuals vs. Indices", ylab = "Residuals", xlab = "Index", col = "blue")

# Check Linearity Assumptions
avPlots(M_selected)
```

**Appendix D (some results of R Code):**

**Resulted models**

# Forward selection model
m1 <- lm(e3_bw ~ h_cohort + hs_c_height_None + e3_sex_None + hs_mo_c_Log2 +
  h_Benzene_Log + h_folic_t1_None + hs_pcb118_cadj_Log2 + hs_cs_m_Log2 +
  hs_dmtp_madj_Log2 + hs_mibp_madj_Log2 + e3_asmokcigd_p_None +
  hs_mbzp_cadj_Log2 + h_builtdens300_preg_Sqrt + hs_no2_dy_hs_h_Log +
  hs_mepa_madj_Log2 + hs_dep_madj_Log2 + hs_pfoa_c_Log2 + hs_c_weight_None +
  hs_pfoa_m_Log2 + hs_dde_madj_Log2 + hs_total_sweets_Ter +
  hs_pbde153_madj_Log2 + h_NO2_Log + hs_pet_dog_r2_None + hs_mnbp_madj_Log2 +
  hs_fastfood_Ter + h_accesspoints300_preg_Log + hs_accesspoints300_h_Log +
  hs_hg_c_Log2 + hs_org_food_Ter + h_temperature_preg_None +
  h_Absorbance_Log + hs_blueyn300_h_None + hs_pm25abs_wk_hs_h_Log +
  hs_hg_m_Log2 + hs_bupa_cadj_Log2 + hs_accesslines300_s_dic0 +
  hs_accesslines300_h_dic0 + hs_popdens_s_Sqrt + hs_greenyn300_s_None +
  hs_hcb_madj_Log2 + hs_mep_madj_Log2 + h_walkability_mean_preg_None +
  h_connind300_preg_Sqrt + hs_connind300_s_Log + hs_dde_cadj_Log2 +
  hs_hcb_cadj_Log2 + hs_tl_cdich_None + hs_pcb153_madj_Log2 +
  hs_pcb170_madj_Log2 + h_PM_Log + hs_dif_hours_total_None +
  hs_cd_m_Log2 + hs_pb_m_Log2 + h_lden_cat_preg_None +
  h_landuseshan300_preg_None +
  hs_landuseshan300_s_None + hs_trafload_h_pow1over3 + hs_fdensity300_h_Log,
  data=predict_data)

# Backward selection model
m2 <- lm(e3_bw ~ hs_c_height_None + hs_c_weight_None + hs_no2_dy_hs_h_Log +
  hs_pm25abs_wk_hs_h_Log + h_accesspoints300_preg_Log + h_builtdens300_preg_Sqrt
  + h_connind300_preg_Sqrt + h_landuseshan300_preg_None +
  h_walkability_mean_preg_None +
  hs_accesspoints300_h_Log + hs_fdensity300_h_Log + hs_connind300_s_Log +
  hs_popdens_s_Sqrt + h_Absorbance_Log + h_Benzene_Log + h_NO2_Log +
  h_PM_Log + hs_KIDMED_None + hs_mvpa_prd_alt_None + hs_dif_hours_total_None +
  hs_cs_m_Log2 + hs_hg_c_Log2 + hs_hg_m_Log2 + hs_mn_c_Log2 +
  hs_mn_m_Log2 + hs_mo_c_Log2 + h_pressure_preg_None + h_temperature_preg_None+
  hs_dde_cadj_Log2 + hs_dde_madj_Log2 + hs_hcb_cadj_Log2 +
  hs_hcb_madj_Log2 + hs_pcb118_cadj_Log2 + hs_pcb153_madj_Log2 +
  hs_pcb170_madj_Log2 + hs_sumPCBs5_madj_Log2 + hs_dep_madj_Log2 +
  hs_dmtp_madj_Log2 + hs_pbde153_madj_Log2 + hs_pfoa_c_Log2 +
  hs_pfoa_m_Log2 + hs_bupa_cadj_Log2 + hs_mbzp_cadj_Log2 +
  hs_mecpp_cadj_Log2 + hs_mep_madj_Log2 + hs_mibp_madj_Log2 +
  hs_mnbp_madj_Log2 + hs_sumDEHP_cadj_Log2 + e3_asmokcigd_p_None +
  h_distinvnear1_preg_Log + hs_trafload_h_pow1over3 + h_cohort +

e3_sex_None + e3_yearbir_None + hs_accesslines300_s_dic0 +
   h_folic_t1_None + h_pavig_t3_None + hs_fastfood_Ter + hs_org_food_Ter +
   hs_pet_dog_r2_None + hs_proc_meat_Ter + hs_total_fruits_Ter +
   hs_total_lipids_Ter + hs_total_sweets_Ter + hs_tl_cdich_None +
   hs_blueyn300_h_None + hs_ln_cat_h_None + hs_lden_cat_s_None +
   hs_globalexp2_None, data=predict_data)

# Stepwise Selection Model
lm(formula = e3_bw ~ hs_c_height_None + hs_c_weight_None + hs_no2_dy_hs_h_Log +
   hs_pm25abs_wk_hs_h_Log + h_accesspoints300_preg_Log + h_builtdens300_preg_Sqrt
   + h_connind300_preg_Sqrt + h_landuseshan300_preg_None +
   h_walkability_mean_preg_None +
   hs_accesspoints300_h_Log + hs_fdensity300_h_Log + hs_connind300_s_Log +
   hs_popdens_s_Sqrt + h_Absorbance_Log + h_Benzene_Log + h_NO2_Log +
   h_PM_Log + hs_KIDMED_None + hs_mvpa_prd_alt_None + hs_dif_hours_total_None +
   hs_cs_m_Log2 + hs_cu_c_Log2 + hs_hg_c_Log2 + hs_hg_m_Log2 +
   hs_mn_c_Log2 + hs_mn_m_Log2 + hs_mo_c_Log2 + hs_pb_m_Log2 +
   h_pressure_preg_None + h_temperature_preg_None + hs_dde_cadj_Log2 +
   hs_dde_madj_Log2 + hs_hcb_cadj_Log2 + hs_hcb_madj_Log2 +
   hs_pcb118_cadj_Log2 + hs_pcb153_madj_Log2 + hs_pcb170_madj_Log2 +
   hs_sumPCBs5_madj_Log2 + hs_dep_madj_Log2 + hs_dmtp_madj_Log2 +
   hs_pbde153_madj_Log2 + hs_pfoa_c_Log2 + hs_pfoa_m_Log2 +
   hs_bupa_cadj_Log2 + hs_mbzp_cadj_Log2 + hs_mecpp_cadj_Log2 +
   hs_mep_madj_Log2 + hs_mibp_madj_Log2 + hs_mnbp_madj_Log2 +
   hs_sumDEHP_cadj_Log2 + e3_asmokcigd_p_None + h_distinvnear1_preg_Log +
   hs_trafload_h_pow1over3 + h_cohort + e3_sex_None + e3_yearbir_None +
   hs_accesslines300_s_dic0 + h_folic_t1_None + h_pavig_t3_None +
   hs_fastfood_Ter + hs_org_food_Ter + hs_pet_dog_r2_None +
   hs_proc_meat_Ter + hs_total_fruits_Ter + hs_total_lipids_Ter +
   hs_total_sweets_Ter + hs_tl_cdich_None + hs_greenyn300_s_None +
   hs_blueyn300_h_None + hs_ln_cat_h_None + hs_lden_cat_s_None +
   Hs_globalexp2_None + hs_tm_mt_hs_h_None, data = predict_data)

# The remaining Added-Value Plots:

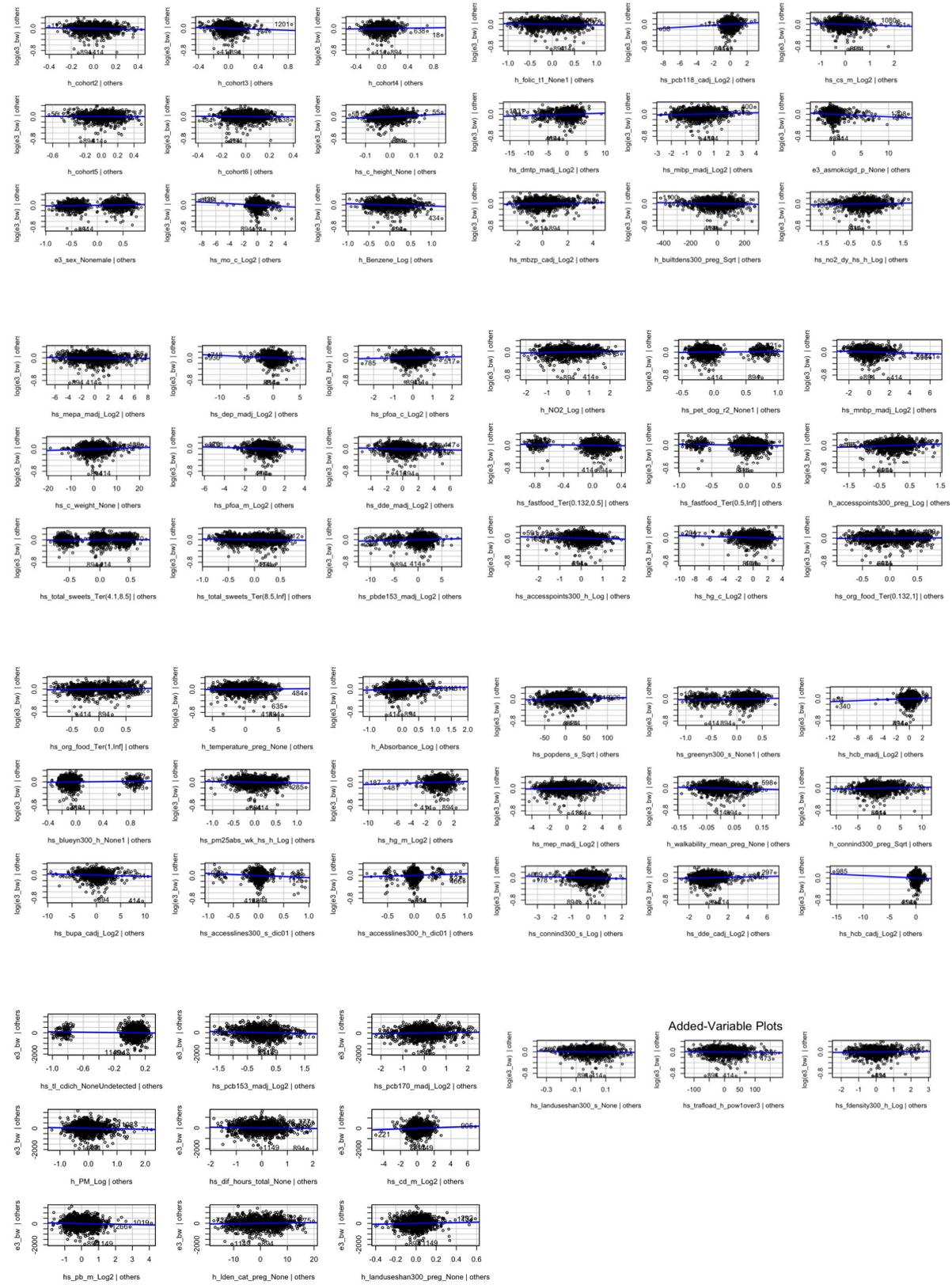